\documentclass[a4paper,11pt,twoside]{article}
\usepackage[dvips]{graphics}
\usepackage[english]{babel}
\usepackage{epsf}
\usepackage{verbatim}
\usepackage{stmaryrd}
\usepackage{latexsym}
\usepackage{makeidx}
\usepackage{epsfig}
\usepackage{amsfonts}
\usepackage{amsmath}
\usepackage{amssymb}
\input xy
\xyoption{all}
\makeindex
\newtheorem{Def}{Definition}[section]
\newtheorem{Prop}[Def]{Proposition}
\newtheorem{Teo}[Def]{Theorem}

\newtheorem{Lem}[Def]{Lemma}

\def\endproof{\ \hfill\hbox{\vbox{\hrule\hbox{\vrule
height5pt\kern5pt\vrule height5pt}\hrule}}\par\medskip\rm}
\title{
\textbf {On $C^0$-variational solutions for Hamilton-Jacobi equations}}

\author{ \sc Olga Bernardi \qquad Franco Cardin \\ { \sc 
}
\\ Dipartimento di Matematica Pura ed Applicata\\ Via Trieste, 63 - 35121 Padova, Italy \\ obern@math.unipd.it \ cardin@math.unipd.it}

\begin{document}

\maketitle

\begin{abstract}
\par \noindent For evolutive Hamilton-Jacobi equations, we propose a
refined definition of $C^0$-variational solution, adapted to Cauchy
problems for continuous initial data. In this weaker framework we
investigate the Markovian (or semigroup) property for these solutions. In
the case of $p$-convex Hamiltonians, when variational solutions are known
to be identical to viscosity solutions, we verify directly the
Markovian property by using minmax techniques. In the non-convex case, we construct an explicit evolutive example where minmax and viscous solutions are different. Provided the initial data allow for the separation of
variables, we also detect the Markovian property for convex-concave Hamiltonians. In this case, and for general initial data, 
we finally give upper and lower Hopf-type estimates for the variational
solutions.
\end{abstract}

\section{Introduction}

The two-fold aim of this paper is to provide a refined definition of variational solution --sometimes called also minmax
solution-- to Cauchy Problems for Hamilton-Jacobi equations of the
evolutive type,
\begin{equation*} (CP)
\begin{cases} \frac{\partial u}{\partial t}(t,x) + H(t,x,\frac{\partial
u}{\partial x}(t,x)) = 0 \\ \\
u(0,x) = \sigma(x),
\end{cases}
\end{equation*}
for the case of continuous initial data $\sigma\in C^0$, and to discuss
some aspects related to the Markovian property of these solutions. The
problem is posed on a smooth, connected and closed manifold $N$ of
dimension $k$ (typically the flat torus $\mathbb{T}^k$), and the
Hamiltonian $H$ is assumed to be of class $C^{1,1}([0,T] \times T^*N)$.

\indent The early notion of a variational solution to $(CP)$ has been
introduced in the Nineties by Chaperon and Sikorav (\cite{CH1}, \cite{V1}
and \cite{VIO}). The construction --see Section \ref{Sezione1}-- employs a
Lusternik-Schnirelman-type procedure based on generating functions
quadratic at infinity (GFQI) for the Lagrangian submanifold --geometric
solution-- $L$ obtained by gluing together the
characteristics of the Hamiltonian vector field for
$\mathcal{H}(t,x,\tau,p) = \tau + H(t,x,p)$ starting from the initial data submanifold
$$\Gamma_\sigma = \{(0,x,-H(0,x,d\sigma(x)),d\sigma(x)): \ x \in N\}
\subset \mathcal{H}^{-1}(0) \subset T^*(\mathbb{R} \times N).$$
In this smooth geometrical environment, the notion of variational solution
--which results to be Lipschitz on finite time, see for example \cite{CH1}
and \cite{PPS}-- overcomes the difficulties arising from the obstruction
to existence of global solutions. 
\par\noindent
However, as the very definition of $\Gamma_{\sigma}$ explains, the global
object $L$ encompassing geometrically the multi-valued features of the
Hamilton-Jacobi problem cannot be defined everywhere for
non-differentiable functions $\sigma$. Consequently, the standard
procedure to obtain the related variational solution does not work
anymore. \\
\indent Our first contribution consists in introducing a natural notion of
$C^0$-variational solution whose Definition \ref{DEFWEAK} is based on a continuity argument inspired
by Viterbo --see Theorem \ref{Viterbo}-- and utilizes $C^1$ approximating
sequences of $\sigma$ in the uniform norm. \\
\indent Extending the notion of variational solution to continuous initial
data is crucial for the investigation of the Markovian (or semigroup)
property, see \cite{VIO}, which means that for all $0 \le
t_1 \le t_2 \le t_3 \le T$ the following holds:
$$J^{t_3,t_2} \circ J^{t_2,t_1} = J^{t_3,t_1},$$
where the map
$$J^{t,t_1}: C^0(N) \to C^0(N), \qquad f(\cdot) \mapsto
J^{t,t_1}(f)(\cdot) := u(t,\cdot),$$
for $t \in [0,T]$ describes the $C^0$-variational solution $u(t,x)$ on the
interval $[0,T]$ starting from $u(t_1,\cdot) = f(\cdot)$. Loosely
speaking, this property requires that the solution at $t_3$ can be
determined by the knowledge of the solution at any intermediate instant
$t_2$.

\indent As is well known, there is another important notion of weak
solution for Hamilton-Jacobi equations with continuous initial data,
namely, the viscosity solution. The reader is referred to \cite{Li}, \cite{Bar}, \cite{BCD} for general review on the theory.
Even thought variational and viscosity solutions have the same analytic
properties, it is not known in the general non-convex case whether they
coincide or not, although they do for $p$-convex Hamiltonians (\cite{JO},
\cite{BeCa}). In Section \ref{Splitting} we construct an evolutive example
showing explicitly the separation between these two notions of solution
for non-convex Hamiltonians.

\indent A main feature of the viscosity solutions is the Markovian
property, see \cite{Li2}. This property has not been proved for
variational solutions in general. By following early suggestions by
Viterbo and Ottolenghi (\cite{VIO}, \cite{V1}), we establish that if the
Hamiltonian is compatible with a natural notion of non-hysteresis for the
above operator $J$, the stronger group, and a fortiori the Markovian,
property is forced on the variational solutions. Specifically, we prove
that the $C^0$-variational solution does not exhibit hysteresis, that is
\begin{equation} \label{ISTE}
J^{t_1,t_2} \circ J^{t_2,t_1} (\sigma) = \sigma \qquad
\text{(non-hysteresis)}
\end{equation}
$\forall t_1,t_2 \in [0,T]$ and $\forall \sigma \in C^0(N)$, if and only
if the group property holds:
\begin{equation} \label{GROUPVERA}
J^{t_3,t_2} \circ J^{t_2,t_1}(\sigma) = J^{t_3,t_1}(\sigma),
\end{equation}
$\forall t_1, t_2, t_3 \in [0,T]$ and $\forall \sigma \in C^0(N)$. \\
\indent We further investigate variational solutions for the case of
$p$-convex Hamiltonians $H(t,x,p) \in C^2([0,T] \times T^*\mathbb{T}^k)$,
when, as recalled, variational and viscous solutions are the same and the
Markovian property holds. However, no direct proof of the latter in a
variational framework was available, and here we prove the Markovian
property for variational solutions by solely using minmax techniques. Our
discussion is useful for generalizations towards some non-convex cases: in
particular, we verify the Markovian property also for convex-concave
Hamiltonians, provided the initial data separate the variables. \\
\indent
Finally, for general initial data, 
we give upper and lower Hopf-type estimates for the variational
solutions. These Hopf-type inequalities are analogous to the
ones given in \cite{BF} and \cite{BO} in the viscous case and they
similarly offer a representation formula for the variational solution
whenever a ``maxmin'' equals to a ``minmax''. \\
\indent This paper is organized as follows. In Sections \ref{Aggiunta} and
\ref{Sezione1} we recall some facts about Lagrangian submanifolds, their
generating functions and the construction of the variational solutions to
$(CP)$. The following section gives a detailed and self-contained proof of
the uniform continuity of the map $J^{t,t_1}$ and the definition of the
$C^0$-variational solution. After a brief comparison with the viscous
case, we also provide an explicit evolutive example where variational and
viscous solutions are actually different. The Markovian property and their
relation to hysteretic phenomena for variational solutions are discussed
in Section \ref{Sezione6}, while in Section \ref{convex} --by utilizing
exclusively minmax arguments-- we prove the Markovian property for
$p$-convex Hamiltonians. In Section \ref{ex} we identify a class of
non-convex problems for which the Markovian property holds and present
upper and lower Hopf-type estimates for variational solutions related to
general convex-concave Hamiltonians. In the Appendix \ref{otto} --essentially following \cite{CH1}-- we explicitly construct the finite-parameters GFQI globally describing the geometric solution $L$ for Hamiltonians $H \in C^{2}([0,T] \times
T^{*}\mathbb{T}^k)$ of the form: $$H(t,x,p) = \frac{1}{2}\langle Ap,p
\rangle + V(t,x,p),$$ $V(t,x,p)$ compactly supported in the $p$ variables,
and show that it is essentially determined by the the quadratic form
$\langle Ap,p \rangle$ involved in the Hamiltonian. As is pointed out in Sections \ref{convex} and \ref{ex}, this structure becomes crucial in the discussion of the Markovian property for variational solutions.
\section{Preliminaries} \label{Aggiunta}
Here we review some topics from the theory of Lagrangian submanifolds and their generating functions. \\
\indent Let $N$ be a smooth, connected and closed (i.e. compact and without boundary) manifold of dimension $k$. We consider the cotangent bundle $T^*N$ equipped with the canonical symplectic form $\omega = dp \wedge dx$, in local coordinates $\omega = \sum_{i=1}^k dp_i \wedge dx^i$. \\
A \textit{Lagrangian submanifold} $L \subset T^{*}N$ is a manifold of dimension $k$ such that $\omega$ vanishes on $L$. A classical argument by Maslov and H\"ormander shows that, at least locally, every Lagrangian submanifold is described by some generating function of the form       
$$S : N \times \mathbb{R}^{h} \longrightarrow \mathbb{R}, \qquad \left(x , \xi\right) \longmapsto S\left(x , \xi\right)$$
in the following way:   
$$L := \left\{ \left(x , \frac{\partial S}{\partial x}\left(x , \xi\right)\right): \ \frac{\partial S}{\partial \xi}\left(x , \xi\right) = 0 \right\},$$
where $0$ is a regular value of the map     
$$\left(x , \xi\right) \longmapsto \frac{\partial S}{\partial \xi}\left(x , \xi\right).$$
We observe that the intersection points of $L$ with the zero section $\mathcal{O}_{T^*N}$ of $T^*N$ are in one-to-one correspondence with the global critical points of $S$. Looking for a condition implying the existence of critical points, the following class of generating functions has been decisive in many issues.
\begin{Def} \label{FGQI} A generating function  $S: N \times \mathbb{R}^h \rightarrow \mathbb{R}$ is quadratic at infinity (GFQI) if for every $\vert \xi \vert > C > 0$
\begin{equation} \label{eq:cond1}
S\left(x , \xi\right) = Q(\xi),
\end{equation}
\par \noindent where $Q(\xi)$ is a nondegenerate quadratic form.
\end{Def}
\par \noindent There are three main operations on generating functions which preserve the quadraticity at infinity property and leave invariant the corresponding Lagrangian submanifolds (see e.g. \cite{We}, \cite{LM}): \\ 
$\bullet$ \textit{Fibered diffeomorphism.}
Let $S : N \times \mathbb{R}^{h} \rightarrow \mathbb{R}$ be a GFQI and $N \times \mathbb{R}^{h} \ni \left(x , \xi\right) \mapsto \left(x , \phi\left(x,\xi\right)\right) \in  N \times \mathbb{R}^{h}$ a map such that, $\forall x \in N$,
$$\mathbb{R}^{k} \ni \xi \longmapsto \phi\left(x , \xi\right) \in \mathbb{R}^{h} $$
\par \noindent is a diffeomorphism. Then
$$S_{1}\left(x , \xi\right) := S\left(x , \phi\left(x,\xi\right)\right)$$
is quadratic at infinity and generates the same Lagrangian submanifold of $S$. \\
$\bullet$ \textit{Stabilization.} 
Let $S : N \times  \mathbb{R}^{h} \rightarrow \mathbb{R}$ be a GFQI Then  
$$S_{1} \left(x , \xi , \eta\right) := S\left(x , \xi\right) + B(\eta),$$ 
where $\eta \in \mathbb{R}^{l}$ and $B(\eta)$ is a nondegenerate quadratic form, generates the same  Lagrangian submanifold of $S$. \\
$\bullet$ \textit{Addition of a constant.} Finally, as a third --although trivial-- invariant operation, we observe that by adding to a generating function $S$ an arbitrary constant $c \in \mathbb{R}$ the described Lagrangian submanifold is invariant. \\ \indent
Crucial problems in the global theory of Lagrangian submanifolds and their parameterizations are $i)$ the existence of a GFQI for a Lagrangian submanifold $L \subset T^{*}N$, $ii)$ the uniqueness of it (up to the three operations described above). \\
The following theorem --see \cite{Si} and \cite{Si1}-- answers to the first question.
\begin{Teo} (Chaperon-Chekanov-Laudenbach-Sikorav) \label{LauSi} 
Let $\mathcal{O}_{T^*N}$ be the zero section of $T^{*}N$ and $\left(\phi_{t}\right)_{t \in [0,1]}$ a Hamiltonian isotopy. Then the Lagrangian submanifold $\phi_{1}\left(\mathcal{O}_{T^*N}\right)$ admits a GFQI
\end{Teo}
\par \noindent The answer to the second problem is due to Viterbo \cite{V2}:
\begin{Teo} (Viterbo) \label{Uniqueness} Let $\mathcal{O}_{T^*N}$ be the zero section of $T^{*}N$ and $\left(\phi_{t}\right)_{t \in [0,1]}$ a Hamiltonian isotopy. Then  the Lagrangian submanifold $\phi_{1}\left(\mathcal{O}_{T^*N}\right)$ admits a unique --that is, up to the above operations-- GFQI
\end{Teo}
Previous theorems --see \cite{Th} and \cite{ThT}-- remain true when the isotopy is only symplectic and we consider generating forms instead of functions. \\
\indent We resume now the Lusternik-Schnirelman calculus of critical values for a GFQI $S(x;\xi)$. Let us consider the sublevel sets
$$S^c := \left\{ (x;\xi) : \ S\left(x;\xi\right) \le c \right\}, \qquad Q^c := \left\{ \xi : \ Q\left(\xi\right) \le c \right\}.$$
We denote by $E^{\pm}$ the positive and negative eigenspaces of $Q$, $D^{\pm}$ large discs in $E^{\pm}$ and $\partial D^{\pm}$ their boundary. \\
Since for $c > 0$ large enough $S^{\pm c} = N \times Q^{\pm c}$, we have
$$H^*(S^c,S^{-c}) = H^*\left(N \times Q^c,N \times Q^{-c}\right)$$
\begin{equation} \label{quadrica}
= H^*(N) \varotimes H^*(Q^c,Q^{-c}) = H^*(N) \varotimes H^*(D^-,\partial D^-).
\end{equation}
Thus, if $\rho$ is the generator of $H^*(D^-,\partial D^-)$, to each cohomology class $\alpha \in H^*(N)$ corresponds the image $T\alpha \in H^*(S^c,S^{-c})$ by the Thom isomorphism:
$$H^*(N) \ni \alpha \mapsto T\alpha = \alpha \wedge \rho \in H^*(S^c,S^{-c}).$$
Let now 
$$i^*_{\lambda}: H^*(S^c,S^{-c}) \to H^*(S^{\lambda},S^{-c})$$
be the pull-back of the inclusion $i_{\lambda}: S^{\lambda} \hookrightarrow S^c$. \\
The idea of utilizing forms in order to construct critical values of $S$ comes back to Birkhoff and Morse. 
\begin{Teo} (Minmax theorem) For $0 \ne \alpha \in H^*(N)$, the minmax value $c(\alpha,S)$ defined as
\begin{equation} \label{alfa}
c(\alpha,S) := \inf \{\lambda \in \left[-c,+c\right] : \ i^*_{\lambda} T\alpha \ne 0\},
\end{equation}
is a critical value of $S$.
\end{Teo}
It is now proved in \cite{VI} that --up to a global shift-- $c(\alpha,S)$ depends only on $L$, not on $S$, and it is thus denoted by $c(\alpha,L)$. \\
\indent We finally recall the definition and properties of the symplectic invariant $\tilde{\gamma}$ for Hamiltonian diffeomorphisms. The presentation of this matter is here oriented to introduce the main tools in order to give a self-contained proof of the forthcoming Theorem \ref{Viterbo}. We refer to \cite{VI}, \cite{CarVit}, \cite{HumT} and \cite{Hum} for an exhaustive treatment of the subject. \\
Let $\mathcal{L}$ be the set of Lagrangian submanifolds of $T^*N$ which are Hamiltonian isotopic to $\mathcal{O}_{T^*N}$ and $L_1, L_2 \in \mathcal{L}$ be generated by the GFQI $S_1(x;\xi)$ and $S_2(x;\eta)$ respectively. We denote by $(S_1 \ \sharp \ S_2)(x;\xi,\eta)$ the GFQI
$$(S_1 \ \sharp \ S_2)(x;\xi,\eta) = S_1(x;\xi) + S_2(x;\eta),$$
and we define
\begin{equation} \label{gamma}
\gamma(L_1,L_2) := c\left(\mu,S_1 \ \sharp \ (-S_2)\right) - c(1,S_1 \ \sharp \ S_2),
\end{equation}
where $1 \in H^0(N)$ and $\mu \in H^k(N)$ are generators. We remind that the value $\gamma(L_1,L_2)$ is often denoted, with abuse of language, by $\gamma(L_1-L_2)$. \\
The previous definition of $\gamma(L_1,L_2)$ pushes to consider the following symplectic invariant for Hamiltonian diffeomorphisms.
\begin{Def} Let $(\phi_t)_{t\in[0,1]}$ be a Hamiltonian isotopy, $\phi = \phi_1$. We set
\begin{equation} \label{gamma-tilde}
\tilde{\gamma}(\phi) := \sup\{\gamma(\phi(L), L): L \in \mathcal{L}\}.
\end{equation}
\end{Def}
All the Hamiltonians are now assumed to be compactly supported, as in \cite{VI}. We refer to \cite{CarVit} for the proof of the following results.
\begin{Prop} \label{VITE} For the function $\tilde{\gamma}$ the following properties hold:
\begin{enumerate}
\item $\tilde{\gamma}(\phi) \ge 0$ and $\tilde{\gamma}(\phi) = 0$ if and only if $\phi = \text{id}$,
\item $\tilde{\gamma}(\phi) = \tilde{\gamma}(\phi^{-1})$,
\item $\tilde{\gamma}(\phi \circ \psi) \le \tilde{\gamma}(\phi) + \tilde{\gamma}(\psi)$ (triangle inequality),
\item $\tilde{\gamma}(\psi \circ \phi \circ \psi^{-1}) = \tilde{\gamma}(\phi)$ (invariance by conjugation).
\end{enumerate}
In particular, $d(\phi_1,\phi_2) := \tilde{\gamma}(\phi_2^{-1} \circ \phi_1)$ defines a metric on the group of Hamiltonian diffeomorphisms of $T^*N$. 
\end{Prop}
\begin{Prop} \label{stima} Assume that $\phi$ is the time-one map associated to the Hamiltonian $H(t,x,p)$. Then we have
$$\tilde{\gamma}(\phi) \le \|H\|_{C^0}.$$
\end{Prop}
We notice that the uniform norm of the Hamiltonian is defined up to the addition of a constant, that is:
$$\|H\|_{C^0} := \sup_{[0,T] \times T^*N} H(t,x,p) - \inf_{[0,T] \times T^*N} H(t,x,p).$$

\section{Construction of the variational solution} \label{Sezione1}
In the sequel, the above topics based on a smooth, connected and closed manifold $N$ will be utilized for a space-time manifold $[0,T] \times N$ of dimension $k+1$. \\
\indent Let consider the following Cauchy problem, related to an evolutive Hamilton-Jacobi equation:
\begin{equation} \label{CA-PR}
\begin{cases} \frac{\partial u}{\partial t}(t,x) + H(t,x,\frac{\partial u}{\partial x}(t,x)) = 0 \\ \\
u(0,x) = \sigma(x)
\end{cases}
\end{equation}
We suppose the Hamiltonian $H : \mathbb{R} \times T^*N \rightarrow \mathbb{R}$ of class $\mathnormal{C}^{1,1}$ and the initial condition $\sigma : N \rightarrow \mathbb{R}$ of class $\mathnormal{C}^1$. \\
Let $\mathbb{R} \times N$ be the ``space-time'', $T^*\left(\mathbb{R} \times N\right)$ its cotangent bundle (endowed with the standard symplectic form $dp \wedge dx + d \tau \wedge dt$):
$$T^*\left(\mathbb{R} \times N\right) = \left\{\left(t,x,\tau,p\right)\right\} \stackrel{\pi_{\mathbb{R} \times N}}{\longrightarrow} \mathbb{R} \times N = \{(t,x)\},$$ 
and $\mathcal{H}\left(t,x,\tau,p\right) = \tau + H\left(t,x,p\right)$. \\
\indent We start resuming the standard procedure to obtain the so-called \textit{geometric solution} to the Cauchy problem (\ref{CA-PR}), which is a Lagrangian submanifold $L \subset T^*N$ attaining --in the sense explained below-- the initial condition $\sigma$ and satisfying: 
$$L \subset \mathcal{H}^{-1}\left(0\right).$$
Let $\Phi_t : \mathbb{R} \times T^*\left(\mathbb{R} \times N\right) \rightarrow T^*\left(\mathbb{R} \times N\right)$ be the flow generated by the Hamiltonian $\mathcal{H} : T^*\left(\mathbb{R} \times N\right) \rightarrow \mathbb{R}$ and $\Gamma_{\sigma}$ the initial data submanifold:
\begin{equation} \label{InitialData}
\Gamma_{\sigma} := \left\{\left(0,x,-H\left(0,x,d\sigma\left(x\right)\right),d\sigma\left(x\right)\right) : \ x \in N\right\} \subset \mathcal{H}^{-1}\left(0\right) \subset T^*\left(\mathbb{R} \times N\right).
\end{equation} 
\par \noindent We note that $\Gamma_{\sigma}$ is the intersection of the Lagrangian submanifold $\Lambda_{\sigma} = \left\{\left(0,x,\tau,d\sigma\left(x\right)\right): \ x \in N, \ \tau \in \mathbb{R}\right\}$ with the hypersurface $\mathcal{H}^{-1}\left(0\right)$:
\begin{equation} \label{int}
\Gamma_{\sigma} = \Lambda_{\sigma} \cap \mathcal{H}^{-1}\left(0\right).
\end{equation}
\begin{Def} \label{GEOMETRIC}
The geometric solution to $(\ref{CA-PR})$ is the submanifold $L = L_{\sigma}$:
\begin{equation} \label{GeoSol}
L := \bigcup_{0 \le t \le T} \Phi_t\left(\Gamma_{\sigma}\right) \subset T^*\left(\mathbb{R} \times N\right).
\end{equation}
\end{Def}
For the following result we refer to \cite{V2}, \cite{CAP} and \cite{BeCa}.
\begin{Prop} \label{SoluzioneGeometrica}
The geometric solution $L$ is an exact Lagrangian submanifold, contained into the hypersurface $\mathcal{H}^{-1}\left(0\right)$ and Hamiltonian isotopic to the zero section $\mathcal{O}_{T^*\left(\left[0,T\right] \times N\right)}$ of $T^*\left(\left[0,T\right] \times N\right)$.
\end{Prop}
As a consequence of the previous Proposition \ref{SoluzioneGeometrica}, Theorem \ref{Uniqueness} guarantees that the Lagrangian submanifold $L$ admits essentially (that is, up to the three operations described above) an unique GFQI 
$$S: (\left[0,T\right] \times N) \times \mathbb{R}^h \rightarrow \mathbb{R}, \qquad \left(t,x;\xi\right) \mapsto S\left(t,x;\xi\right).$$ 
Moreover, we can assume that the graph of $S\left(t,x;\xi\right)$ at $t = 0$ coincides with $\Gamma_{\sigma}$:
$$\Gamma_{\sigma} = L \cap \pi^{-1}_{[0,T] \times N} \left(\{0\} \times N\right)$$
$$= \left\{ \left(0,x,\frac{\partial S}{\partial t}\left(0,x;\xi\right),\frac{\partial S}{\partial x}\left(0,x;\xi\right)\right) : \ \frac{\partial S}{\partial \xi}\left(0,x;\xi\right) = 0 \right\}.$$
The quadraticity at infinity property is crucial: as explained here below, variational solutions arise from the application of the Lusternik-Schnirelman method to $S\left(t,x;\xi\right)$. \\
\indent We denote by $S_{(t,x)}(\xi) = S(t,x;\xi)$ the restriction of $S$ to the fiber over $(t,x)$ and we look for a minmax value of the function $S_{(t,x)}$. Since the cohomology of the point is one dimensional, denoting by $1_{(t,x)}$ its generator, we give the following
\begin{Def} \label{Minimax soluzione} The variational solution of $\left(\ref{CA-PR}\right)$ is the function --see (\ref{alfa})-- 
\begin{equation} \label{MiniMax}
\left(t,x\right) \mapsto u(t,x) := c(1_{(t,x)},S_{(t,x)}).
\end{equation}
\end{Def}
\par \noindent The next fundamental theorem has been proved by Chaperon in \cite{CH1}, see also \cite{VIO}. A simple and self-contained proof of the Lipschitzianity of $u(t,x)$ is also given in \cite{PPS}.
\begin{Teo} \label{CON-LIP}
The variational solution $u(t,x)$ is a weak solution to $\left(CP\right)$, Lipschitz on finite time, which does not depend on the choice of the GFQI $S(t,x;\xi)$.
\end{Teo}
We observe that the definition of variational solution arises naturally in the compact case, when the Theorem \ref{Uniqueness} is satisfied. Moreover, since $H^*(D^-,\partial D^-) = H^*(Q^c,Q^{-c})$ (see (\ref{quadrica})), for any point $(t,x) \in [0,T] \times N$, the minmax critical value $u(t,x)$ is determined by the Morse index of the quadratic form $Q(\xi)$ coinciding with $S(t,x;\xi)$ out of a compact set in the parameters $\xi$. As it is largely known, for $p$-convex Hamiltonians --we refer to \cite{JO} and \cite{BeCa} for details-- the Morse index of the quadratic form $Q(\xi)$ is $0$. Then
\begin{equation} \label{min}
u(t,x) = \min_{\xi \in \mathbb{R}^h} S(t,x;\xi).
\end{equation}
Analogously, recalling that $c(\mu,S) = -c(1,-S)$ (see \cite{VI}), the $p$-concave case gives:
\begin{equation} \label{max}
u(t,x) = -\min_{\xi \in \mathbb{R}^h} -S(t,x;\xi) = \max_{\xi \in \mathbb{R}^h} S(t,x;\xi).
\end{equation}
A rather general discussion for Hamiltonians of the form $H(t,x,p) = \frac{1}{2}\langle Ap,p \rangle + V(t,x,p)$ is performed in the Appendix. The previous representation formulas (\ref{min}) and (\ref{max}) will be crucial in Section \ref{sette}, which is devoted to the discussion of the Markovian property for variational solutions in the convex and convex-concave cases. \\
We conclude with the following result, see \cite{HumT}, showing the $\gamma$-continuity of the variational solution with respect to the geometric one.
\begin{Prop} \label{HU} Let $L_1, L_2$ and $u_1, u_2$ be the geometric and variational solutions for the Cauchy problems referred to the initial data $\sigma_1$ and $\sigma_2$ respectively. Then we have
$$\|u_1 - u_2\|_{C^0} \le \gamma(L_1,L_2).$$
\end{Prop}

\section{$C^0$-variational solutions and viscosity solutions} \label{Sezione5}
In the previous section, given a pair $(H,\sigma) \in C^{1,1}([0,T] \times T^*N) \times C^1(N)$, we have outlined the construction of the geometric solution $L$ (see (\ref{GeoSol})) and the variational solution $u$ (see (\ref{MiniMax})), which results (Theorem \ref{CON-LIP}) Lipschitz on finite time. Thus we consider the application $J$:
$$J:  C^{1,1}([0,T] \times T^*N) \times C^1(N) \rightarrow C^{0,1}([0,T] \times N)$$
\begin{equation} \label{J}
(H,\sigma) \mapsto u =: J(H,\sigma)
\end{equation}
associating to the pair $(H,\sigma)$ the variational solution to the Cauchy problem (\ref{CA-PR}). Often in the following we refer to the above solution as to the $C^{0,1}$-variational solution. \\
The very construction of $u$ enables us to establish the next main result outlined by Viterbo in \cite{VIO} and \cite{V2}. Here below, we focus on a detailed and self-contained proof of the uniform continuity of $J$ with respect to the initial data $\sigma$.
\begin{Teo} \label{Viterbo}
The application $J$ is uniformly continuous if all the spaces are equipped with the $C^0$ topology. Thus it extends to an uniformly continuous map, still denoted by J
$$J: C^{0,1}([0,T] \times T^*N) \times C^0(N) \rightarrow C^0([0,T] \times N).$$
\end{Teo} 
\textit{Proof.} In \cite{VIO} is proved that: 
$$\|J(H_1,\sigma) - J(H_2,\sigma)\|_{C^0} \le T \|H_1 - H_2\|_{C^0}$$
and the extension of $J$ to Lipschitz Hamiltonians is largely discussed. Here we investigate on the following estimate involving the initial data:
$$\|J(H,\sigma_1) - J(H,\sigma_2)\|_{C^0} \le \|\sigma_1 - \sigma_2\|_{C^0}.$$
For $i = 1,2$, let $\Gamma_{\sigma_i}$, $L_{\sigma_i}$ and $S_i$ be respectively the initial data submanifold, the geometric solution and the corresponding GFQI related to the initial datum $\sigma_i$. Moreover, since the geometric solution --see Proposition \ref{SoluzioneGeometrica}-- is Hamiltonian isotopic to the zero section $\mathcal{O} = \mathcal{O}_{T^*([0,T] \times N)}$ of $T^*([0,T] \times N)$, we denote by $\rho_1$ and $\rho_2$ the time-one map such that
\begin{equation} \label{LO}
L_{\sigma_1} =  \bigcup_{0 \le t \le T} \Phi_t\left(\Gamma_{\sigma_1}\right) = \rho_1(\mathcal{O}), \qquad L_{\sigma_2} =  \bigcup_{0 \le t \le T} \Phi_t\left(\Gamma_{\sigma_2}\right) = \rho_2(\mathcal{O}).
\end{equation}
Hence, as a consequence of Proposition \ref{HU}, we have:
$$\|J(H,\sigma_1) - J(H,\sigma_2)\|_{C^0} \le \gamma(L_{\sigma_1},L_{\sigma_2}) = \gamma(\rho_1 \circ \rho_2^{-1}(L_{\sigma_2}), L_{\sigma_2}) \le \tilde{\gamma}(\rho_1 \circ \rho_2^{-1}).$$
In order to estimate $\tilde{\gamma}(\rho_1 \circ \rho_2^{-1})$, we proceed into two steps. We first explicitly construct a Hamiltonian isotopy $(\Psi_s)_{s \in [0,1]}$ of $T^{*}(\mathbb{R} \times N)$ with the following two properties:
$$\Psi_1(\Gamma_{\sigma_1}) = \Gamma_{\sigma_2}, \qquad \tilde{\gamma}(\Psi_1) \le \|\sigma_1 - \sigma_2 \|_{C^0}.$$ 
Secondly, we prove that $\tilde{\gamma}(\rho_1 \circ \rho_2^{-1}) = \tilde{\gamma}(\Psi_1)$. \\
For $i = 1,2$, let us consider the submanifold:
$$\Gamma_{\sigma_i}^0 = \Gamma_{\sigma_i} \cap \{\tau=0\} = \{(0,x,0,d\sigma_i(x)): \ x \in N\}.$$ 
In order to move $\Gamma_{\sigma_1}$ into $\Gamma_{\sigma_2}$, we start observing that $\Gamma_{\sigma_1}^0$ can be easily moved into $\Gamma_{\sigma_2}^0$ by using the time-one flow of the Hamiltonian: 
$$K(t,x,\tau,p) = \sigma_1(x) - \sigma_2(x).$$ 
Therefore, it only remains to move $\Gamma_{\sigma_1}$ into $\Gamma_{\sigma_1}^0$ and $\Gamma_{\sigma_2}^0$ into $\Gamma_{\sigma_2}$ and this can be achieved by employing $C^0$ arbitrary small Hamiltonians. In particular, given a function $g(t) \in C^{1}([0,T])$ with $g(0) = 0$ and $g'(0) = 1$ (we will choose $\| g \|_{C^0} <<1$ below), the time-one flow of the Hamiltonian: 
$$K_1(t,x,\tau,p) = -g(t) H(0,x,d\sigma_1(x))$$ 
moves $\Gamma_{\sigma_1}$ into $\Gamma_{\sigma_1}^0$ and the time-one flow of the Hamiltonian:
$$K_2(t,x,\tau,p) := g(t) H(0,x,d\sigma_2(x))$$ 
moves $\Gamma_{\sigma_2}^0$ into $\Gamma_{\sigma_2}$. Hence, the required Hamiltonian isotopy $(\Psi_s)_{s \in [0,1]}$ is the composition of the three Hamiltonian flows constructed just above. Using now the triangle inequality for the function $\tilde{\gamma}$ (see Proposition \ref{VITE}) and the estimate given in Proposition \ref{stima} and since $\|g\|_{C^0}$ is arbitrary small, we conclude that $\tilde{\gamma}(\Psi_1) \le \|\sigma_1 - \sigma_2 \|_{C^0}$. \\
Here below we show that $\| \rho_1 \circ \rho_2^{-1} \|_{C^0} = \| \Psi_1 \|_{C^0}$. From the one hand, taking into account (\ref{LO}):
$$L_{\sigma_1} = \rho_1(\mathcal{O}) = \rho_1 \circ \rho_2^{-1}(L_{\sigma_2})$$
$$= \rho_1 \circ \rho_2^{-1} \left(\bigcup_{0 \le t \le T} \Phi_t\left(\Gamma_{\sigma_2}\right)\right) = \bigcup_{0 \le t \le T} \rho_1 \circ \rho_2^{-1} \circ \Phi_t\left(\Gamma_{\sigma_2}\right).$$
From the other hand, since $\Psi_1(\Gamma_{\sigma_1}) = \Gamma_{\sigma_2}$:
$$L_{\sigma_1} = \bigcup_{0 \le t \le T} \Phi_t\left(\Gamma_{\sigma_1}\right) =  \bigcup_{0 \le t \le T} \Phi_t \circ \Psi_1^{-1} \left(\Gamma_{\sigma_2}\right).$$
Therefore, $\forall t \in [0,T]$:
$$\rho_1 \circ \rho_2^{-1} \circ \Phi_t \left(\Gamma_{\sigma_2}\right) = \Phi_t \circ \Psi_1^{-1} \left(\Gamma_{\sigma_2}\right),$$
implying that:
$$\Phi_t^{-1} \circ \rho_1 \circ \rho_2^{-1} \circ \Phi_t \left(\Gamma_{\sigma_2}\right) = \Psi_1^{-1} \left(\Gamma_{\sigma_2}\right).$$
Moreover, from the very construction of $\Psi_1$ we obtain that the same property holds for the Lagrangian submanifold of $T^*([0,T] \times N)$:
$$\Lambda := \{(0,x,e-H(0,x,d\sigma_2(x)),d\sigma_2(x)): \ x \in N, \ e \in \mathbb{R}\},$$
that is
\begin{equation*}
\Phi_t^{-1} \circ \rho_1 \circ \rho_2^{-1} \circ \Phi_t \left(\Lambda\right) = \Psi_1^{-1} \left(\Lambda\right).
\end{equation*} 
Therefore, keeping in mind the invariance by conjugation and inverse of $\tilde{\gamma}$ (see Proposition \ref{VITE}), we conclude that $\forall t \in [0,T]$:
$$\tilde{\gamma}(\Phi_t^{-1} \circ \rho_1 \circ \rho_2^{-1} \circ \Phi_t ) = \tilde{\gamma}(\rho_1 \circ \rho_2^{-1}) = \tilde{\gamma}(\Psi_1^{-1}) = \tilde{\gamma}(\Psi_1) \le \|\sigma_1 - \sigma_2 \|_{C^0},$$
and the required estimate is proved. \hfill $\Box$ \\ \\
\indent Previous theorem allows us to consider Cauchy problems with weakly regular initial data. In fact, for only continuous functions $\sigma$, the above --see (\ref{InitialData})-- initial data submanifold $\Gamma_{\sigma}$ cannot be defined anymore and consequently the standard procedure to obtain the geometric solution $L$ and the related variational solution $u$ does not work. However, every continuous initial datum $\sigma \in C^0(N)$ can be approximated in the uniform convergence by a sequence of differentiable functions $\sigma_n \in C^1(N)$, for which we construct the related variational solution $J(H,\sigma_n) = u_{\sigma_n}$. As a consequence of the continuity of $J$, it is easy to prove by direct computation that $i)$ $u_{\sigma_n}$ is a Cauchy sequence --therefore it converges on the complete space $C^0([0,T] \times N)$-- and $ii)$ its limit does not depend on the chosen approximating sequence $\sigma_n$. These arguments justify the next 
\begin{Def} \label{DEFWEAK} ($C^0$-variational solution) Given a continuous initial datum $\sigma \in C^0(N)$, the $C^0$-variational solution for the Cauchy problem (\ref{CA-PR}) is the unique function $u_\sigma \in C^0([0,T] \times N)$ such that, for any arbitrary $C^1$ approximating sequence $\sigma_n$:
$$C^1(N) \ni \sigma_n \stackrel{C^0}{\longrightarrow} \sigma \in C^0(N),$$
with related $C^{0,1}$-variational solutions $J(H,\sigma_n) = u_{\sigma_n}$, we have that
\begin{equation} \label{limit}
\lim_{n \to +\infty}\| u_{\sigma_n} - u_\sigma \|_{C^0} = 0 \qquad \text{on } [0,T] \times N.
\end{equation}
\end{Def}

\indent Since the beginning of the Eighties of the past century, it is present in literature a rather fruitful weak notion of solution for Hamilton-Jacobi equations with continuous initial data, namely, the \textit{viscosity solution}. We refer to \cite{Li}, \cite{Bar} and \cite{BCD} for a detailed review on the subject. 
\begin{Def} \label{Visco}
A function $u \in C((0,T) \times N)$ is a viscosity subsolution [supersolution] of 
\begin{equation} \label{HJ}
\frac{\partial u}{\partial t}(t,x) + H(t,x,\frac{\partial u}{\partial x}(t,x)) = 0
\end{equation}
if, for any $\phi \in C^1((0,T) \times N)$,
\begin{equation} \label{subsol}
\frac{\partial \phi}{\partial t}(\bar{t},\bar{x}) + H(\bar{t},\bar{x},\frac{\partial \phi}{\partial x}(\bar{t},\bar{x})) \le 0 \qquad [\ge 0]
\end{equation}
\noindent at any local maximum [minimum] point $(\bar{t},\bar{x}) \in (0,T) \times N$ of $u - \phi$. Finally, $u$ is a viscosity solution of (\ref{HJ}) if it is simultaneously a viscosity subsolution and supersolution. 
\end{Def}
The two notions of solution have the same analytic properties, that is, theorems of existence and uniqueness do hold, see \cite{VIO}. It is also remarkable that  the Definition \ref{DEFWEAK} of $C^0$-variational solution is consistent with the well-known $C^0$-stability (for which we refer to \cite{Li}) of viscosity solutions. Variational and viscous solutions coincide when the Hamiltonian is $p$-convex, see \cite{JO}. A detailed proof of this fact in the case $H(x,p) = \frac{1}{2} |p|^2 + V(x), \ x \in\mathbb{T}^k$, is also given in \cite{BeCa} using the Lax-Oleinik representation formula. However, we do not know in literature examples showing the splitting between minmax and viscosity solutions: a task of the next section is to construct an explicit evolutive example marking the separation between the two notions of solution.

\subsection{Splitting between minmax and viscosity solutions: an evolutive example} \label{Splitting}
In the paper \cite{MCC}, starting from the non-convex Hamiltonian 
$$H(x,p) = p - p^3 - x,$$ 
it is explicitly shown that the minmax selection from the GFQI of a Lagrangian submanifold contained in $H^{-1}(0)$ cannot be a subsolution for the corresponding \textit{non-evolutive} equation. The evolution setting for the same Hamiltonian displays some intriguing aspects connected to the choice of the initial data, e.g. the Cauchy problem:
\begin{equation} \label{EXX}
\begin{cases}
\frac{\partial u}{\partial t}(t,x) + \frac{\partial u}{\partial x}(t,x) - \left(  \frac{\partial u}{\partial x}(t,x) \right)^3 - x = 0 \\ \\
u(0,x) = \sigma(x) = 0
\end{cases}
\end{equation}
admits a global classical solution (which is both the variational and the viscous one) and does not display any separation between the two notions of solution. Nevertheless, drawing on the above non-evolutive case, we explicitly construct an \textit{evolutive} example where minmax  and viscous solutions actually differ. \\
\indent Alternatively to the problem (\ref{EXX}), we consider non-vanishing initial impulses $d\sigma(x) = v(x)$ such that, for fixed $\varepsilon > 0$:
$$v(x)  = \begin{cases}
\text{the positive solution of } v - v^3 - x = 0 & \text{ if } x < -\varepsilon \\ \\
\text{the negative solution of } v - v^3 - x = 0 & \text{ if } x > \varepsilon \\ \\
\text{any monotone smooth joint of the above branches} & \text{ if } |x| \le -\varepsilon
\end{cases}$$
and we investigate on the solution $u(t,x)$ for the initial datum $u(0,x) = \sigma(x)$. \\
Hamilton's equations related to $\mathcal{H}(t,x,\tau,p) = \tau + p - p^3 - x$ are
\begin{equation} \label{HAMM}
\begin{cases}
\dot{t} = 1 \\
\dot{x} = 1 - 3p^2 \\
\dot{\tau} = 0 \\
\dot{p} = 1
\end{cases}
\end{equation}
Focusing on $|x_0| > \varepsilon$, the corresponding initial data submanifold is given by
$$\Gamma_{\sigma, |x_0| > \varepsilon} = \left\{\left(0,x_0, -H(x_0,d\sigma(x_0)), d\sigma(x_0) \right): \ |x_0| > \varepsilon \right\}$$
$$= \left\{\left(0,x_0,0,v(x_0)\right): \ |x_0| > \varepsilon \right\}.$$
Since the flow $\phi^t_H$ reads
\begin{equation} \label{evoluti}
\begin{cases} 
x(t) = x_0 + t - 3v^2(x_0)t - 3v(x_0)t^2 -t^3 \\
p(t) = v(x_0) + t,
\end{cases}
\end{equation}
starting from $ |x_0| > \varepsilon$ and taking into account the very property of the function $v(x_0)$ again, we obtain
$$\begin{cases}
x(t) = v(x_0) + t - (v(x_0) + t)^3 \\
p(t) = v(x_0) + t
\end{cases}$$
so that the corresponding geometric solution equals
\begin{equation} \label{GeoNostra}
L_{|x_0| > \varepsilon} = \bigcup_{0 \le t \le T} \left\{\left(t,v(x_0) + t - (v(x_0) + t)^3,0, v(x_0) +t\right): \ |x_0| > \varepsilon \right\}.
\end{equation}
A simple recognition of $\text{(}\ref{evoluti}\text{)}_1$ shows that the $x$-components of the characteristics starting from $|x_0| \le \varepsilon$ remain definitively away from $0$ for $t>0$ big enough. As a consequence, locally to $x = 0$ and for $t> 0$ sufficiently large, the geometric solution corresponds just to (\ref{GeoNostra}). There, the Lagrangian submanifold is three-valued and it is locally  described by the GFQI 
\begin{equation} \label{S-Nostra}
S(t,x;\xi) = \frac{1}{2}\xi^2 + t\xi - \frac{3}{4}(\xi-x+t)^{\frac{4}{3}} + \frac{1}{2}t^2.
\end{equation}
Due to the presence of the leading term $\frac{1}{2}\xi^2$, the minmax procedure for $S(t,x;\xi)$ selects the minimum value over the parameters $\xi$. More precisely, from the condition
$$\frac{\partial S}{\partial \xi} (t,x;\xi) = \xi + t - (\xi-x+t)^{\frac{1}{3}} = 0,$$
or equivalently $(\xi + t) - (\xi + t)^3 - x = 0$, we obtain directly that the corresponding minmax solution $u(t,x)$ assumes the following analytic form. \\
For $x < 0$, let us denote by $v^+(x)$ the unique positive root of the cubic polynomial $v - v^3 - x$. Then  we have that:
$$u(t,x) = S(t,x;v^+(x) - t)$$ 
$$= \frac{1}{2}(v^+(x) - t)^2 + t(v^+(x) - t) - \frac{3}{4}(v^+(x) -x)^{\frac{4}{3}} + \frac{1}{2}t^2.$$ 
For $x > 0$, if we indicate by $v^-(x)$ the unique negative root of $v - v^3 - x = 0$, 
$$u(t,x) = S(t,x;v^-(x) - t)$$ 
$$= \frac{1}{2}(v^-(x) - t)^2 + t(v^-(x) - t) - \frac{3}{4}(v^-(x) -x)^{\frac{4}{3}} + \frac{1}{2}t^2.$$
Finally, it can be verified directly that the minimum $u(t,0) = -\frac{1}{4}$ is attained for the values $1 - t$ and $-1-t$ of the parameter $\xi$. From the derivatives
$$\frac{\partial S}{\partial t}(t,x;\xi) = \xi - (\xi-x+t)^{\frac{1}{3}} + t, \qquad \frac{\partial S}{\partial x}(t,x;\xi) = (\xi-x+t)^{\frac{1}{3}},$$
we obtain that the function $u(t,x)$ results non-differentiable in the points $(t,0)$, with corresponding superdifferential $D^+u(t,0) = \{ 0 \} \times [-1,1]$ and subdifferential $D^-u(t,0) = \emptyset$. \\
By standard arguments --see \cite{BCD}-- we conclude that, at any point $(t,0)$ as above, the function $u(t,x)$ cannot be a viscosity solution to the evolutive Hamilton-Jacobi equation $\text{(}\ref{EXX}\text{)}_1$. 
In fact, taking e.g. the value $(0,\frac{1}{\sqrt{3}}) \in D^+u(t,0)$ and since $H(0,\frac{1}{\sqrt{3}}) = \frac{2}{3\sqrt{3}} > 0$, the subsolution property fails.  

\section{The Markovian property} \label{Sezione6}
The enlargement of the notion of variational solution to continuous initial data turns out crucial in order to deal with the \textit{Markovian property} for such solutions, which can be explained as follows. \\
\indent Let $H$ be a fixed $C^{1,1}$ Hamiltonian and $f$ a continuous initial datum. Let denote by $u \in C^{0}([0,T] \times N)$ the variational solution (related to the refined Definition \ref{DEFWEAK}) with datum $f$ at time $t_1 \in [0,T]$:
\begin{equation} \label{mar}
\begin{cases}
\frac{\partial u}{\partial t}(t,x) + H(t,x,\frac{\partial u}{\partial x}(t,x)) = 0 \\ \\
u(t_1,x) = f(x)
\end{cases}
\end{equation}
Equivalently, the above solution $u(t,x)$ is described by the following map $J^{t,t_1}$: 
\begin{equation} \label{J12}
J^{t,t_1}: C^{0}(N) \rightarrow C^{0}(N), \qquad f(\cdot) \mapsto J^{t,t_1}(f)(\cdot) := u(t,\cdot),
\end{equation}
where $t \in [0,T]$. The Markovian (or semigroup) property is one of the main features of the viscosity solutions; it means that:
\begin{equation} \label{MARKVERA} 
J^{t_3,t_2} \circ J^{t_2,t_1}(\sigma) = J^{t_3,t_1}(\sigma),
\end{equation}
$\forall 0 \le t_1\le t_2 \le t_3 \le T$ and $\forall \sigma \in C^0(N)$. A general treatment of the nonlinear semigroup associated to a first-order Hamilton-Jacobi equation with uniformly continuous initial data is developed in \cite{NL} and \cite{Li2} for viscosity solutions. \\
This is not the case of variational solutions: Viterbo and Ottolenghi in \cite{VIO} (see also \cite{V1}) suggest that the failure of this property can be marked by a sort of  ``hysteresis phenomena''. In the next proposition, we show that  for the variational solutions a natural mathematical notion of non-hysteresis --see (\ref{ISTE})-- is actually forcing the stronger group property.
\begin{Prop} \label{SeSoloSe} Let $H \in C^{1,1}([0,T] \times T^*N)$. The $C^0$-variational solution does not exhibit hysteresis, that is
\begin{equation} \label{ISTE}
J^{t_1,t_2} \circ J^{t_2,t_1} (\sigma) = \sigma \qquad \text{(non-hysteresis)}
\end{equation}
$\forall t_1,t_2 \in [0,T]$ and $\forall \sigma \in C^0(N)$, if and only if the group property holds: 
\begin{equation} \label{GROUPVERA} 
J^{t_3,t_2} \circ J^{t_2,t_1}(\sigma) = J^{t_3,t_1}(\sigma),
\end{equation}
$\forall t_1, t_2, t_3 \in [0,T]$ and $\forall \sigma \in C^0(N)$. 
\end{Prop}
\textit{Proof.} We have to prove only the sufficiency of the condition (\ref{ISTE}), the necessity is in fact immediate. Let $u \in C^0([0,T] \times N)$ be the variational solution to the Cauchy problem:
$$
\begin{cases}
\frac{\partial u}{\partial t}(t,x) + H(t,x,\frac{\partial u}{\partial x}(t,x)) = 0 \\ \\
u(t_1,x) = \sigma(x)
\end{cases}
$$
and $f(\cdot) := u(t_2,\cdot) \in C^0(N)$. We consider an arbitrary $C^1$ approximating sequence $f_n \stackrel{C^0}{\to} f$, generating the $C^{0,1}$-variational solutions $u_{f_n}$ on $[0,T] \times N$ with $u_{f_n}(t_2,\cdot) = f_n(\cdot)$. Under the Definition \ref{DEFWEAK}, the variational solution $u_f \in C^0([0,T] \times N)$ to the Cauchy problem starting from $u_f(t_2, \cdot) = f(\cdot)$ is given by:
\begin{equation} \label{ZERO}
\lim_{n \to +\infty} \|u_{f_n} - u_f\|_{C^0} = 0,
\end{equation}
and $u_f(t_3,\cdot)$ is exactly:
\begin{equation} \label{PRIMA}
u_f(t_3,\cdot) = J^{t_3,t_2} \circ J^{t_2,t_1} (\sigma)(\cdot).
\end{equation}
We proceed by introducing
$$\Sigma_n(\cdot) := J^{t_1,t_2}(f_n)(\cdot) \in C^{0,1}(N),$$
with related $C^{0,1}$-variational solution $u_{\Sigma_n}$. \\
Since $f_n$ is convergent and $J^{t_1,t_2}$ is continuous, the sequence $\Sigma_n$ is convergent and, as a consequence of the non-hysteresis assumption (\ref{ISTE}), $\Sigma_n \stackrel{C^0}{\to} \sigma$ does hold. Moreover, as by construction both $u_{\Sigma_{n}}$ and $u_{f_n}$ produce the solution to the Cauchy problem:
$$
\begin{cases}
\frac{\partial u}{\partial t}(t,x) + H(t,x,\frac{\partial u}{\partial x}(t,x)) = 0 \\ \\
u(t_2,x) = f_n(x)
\end{cases}
$$
by uniqueness, we conclude that $u_{\Sigma_{n}} = u_{f_n}$ on $[0,T] \times N$. \\
For every $\Sigma_n \in C^{0,1}(N)$, we take now a $C^1$ regularizing sequence $\Sigma_{n,m} \stackrel{C^0}{\to} \Sigma_n$ and we determine $J^{t_3,t_1}(\sigma)(\cdot)$ by utilizing the (diagonal) sequence: 
$$\sigma_n := \Sigma_{n,n} \stackrel{C^0}{\to} \sigma,$$ 
with corresponding $C^{0,1}$-variational solutions $u_{\sigma_n}$. Under the Definition \ref{DEFWEAK}, the variational solution $u_{\sigma} \in C^0([0,T] \times N)$ is given by:
\begin{equation} \label{NUOVA}
\lim_{n \to +\infty} \|u_{\sigma_{n}} - u_{\sigma}\|_{C^0} = 0,
\end{equation}
and $u_\sigma(t_3,\cdot)$ is exactly:
\begin{equation} \label{TERZA}
u_{\sigma}(t_3,\cdot) = J^{t_3,t_1} (\sigma)(\cdot).
\end{equation}
However, from
$$\| u_{\sigma_n} - u_f \|_{C^0} \le \| u_{\sigma_n} - u_{\Sigma_{n}} \|_{C^0} + \| u_{\Sigma_{n}} - u_{f_n} \|_{C^0} + \| u_{f_n} - u_{f} \|_{C^0},$$
and by using (\ref{ZERO}) and the fact that $u_{\Sigma_{n}} = u_{f_n}$ on $[0,T] \times N$, we achieve that
\begin{equation} \label{QUARTA}
\lim_{n \to +\infty} \|u_{\sigma_n} - u_f\|_{C^0} = 0.
\end{equation}
As a consequence, compare (\ref{ZERO}), (\ref{NUOVA}) and (\ref{QUARTA}), $u_f$ and $u_{\sigma}$ coincide on $[0,T] \times N$ and for $t = t_3$ we obtain (\ref{GROUPVERA}). \hfill $\Box$ 

\section{Markovian variational solutions} \label{sette}
The above weak Definition \ref{DEFWEAK} is substantially based on the construction of converging sequences of minmax solutions:
$$u_{\sigma_n} \stackrel{C^{0}}{\rightarrow} u_{\sigma},$$
where $C^1(\mathbb{T}^k) \ni \sigma_n \stackrel{C^{0}}{\rightarrow} \sigma \in C^0(\mathbb{T}^k)$. We remind that, for a smooth (i.e. at least $C^1$) initial datum, explicit formulas for variational solutions pass through the construction of the finite parameters GFQI globally describing the corresponding geometric solution (see (\ref{GeoSol})): the general procedure is based on a finite reduction of the Hamilton-Helmholtz functional, which can be performed by the ``broken geodesics'' method of Chaperon (see \cite{CH1}, \cite{CH2} and also \cite{Si}) or, alternatively, by an Amann-Conley-Zehnder reduction procedure, see for example \cite{Aeb} and \cite{Car}. \\
In the Appendix --substantially following the line of thought of Chaperon-- we resume such a construction in a rather general case, that is for Hamiltonians $H(t,x,p) \in C^{2}([0,T] \times T^{*}\mathbb{T}^k)$ of the form: 
\begin{equation} \label{HAM}
H(t,x,p) = \frac{1}{2}\langle Ap,p \rangle + V(t,x,p),
\end{equation} 
where $A^t = A$, $\text{det}(A) \ne 0$ and $V(t,x,p)$ is compactly supported in the $p$ variables. \\
\indent We assume the global existence of the Legendre transformation: 
\begin{itemize}
\item[$(\star)$] the map $p \mapsto D_pH(t,x,p)$ provides a global diffeomorphism of $\mathbb{R}^k$ into itself, uniformly Lipschitz with its inverse.
\end{itemize} 
The previous condition is in particular satisfied when $H(t,x,p)$ is $p$-convex:
\begin{equation} \label{con}
\exists c > 0: \qquad \langle D^2_pH(t,x,p) \lambda, \lambda \rangle \ge c|\lambda|^2,
\end{equation}
$\forall \lambda \in \mathbb{R}^k, \ \forall (t,x,p) \in [0,T] \times T^*\mathbb{T}^k$. However, it includes  also a large class of non-convex cases: for instance $H(p) = \frac{1}{2}\langle Ap,p \rangle$ with general hyperbolic matrices $A$. \\
We refer to the Appendix for the proof of the next result.
\begin{Teo} \label{AA} 
Let $H(t,x,p) \in C^{2}([0,T] \times T^{*}\mathbb{T}^k)$ be a Hamiltonian function of the form (\ref{HAM}) and $\sigma \in C^1(\mathbb{T}^k)$. \\
We suppose $H(t,x,p)$ satisfying condition $(\star)$ and we denote by $S^t_0(X,x;U)$, $U \in \mathbb{R}^{Nk}$, the generating function for the flow $\phi^{t,0}_H: (X,P) \mapsto (x,p)$. \\
Then the (broken geodesics) generating function
\begin{equation} \label{ESSE t}
S^t(x;\xi,U) := \sigma(\xi) + S^t_0(\xi,x;U)
\end{equation}
for the Lagrangian wavefront $\phi^{t,0}_H(\text{Im}(d\sigma))$ is quadratic at infinity, with quadratic form given by the nondegenerate $(N+1)k \times (N+1)k$ matrix:
$$\textbf{A} = \begin{pmatrix}
A & 0 & 0 \\
0  & \ddots & 0 \\
0  & 0 & A
\end{pmatrix}.$$
\end{Teo}
As a consequence, in the hypothesis of the previous theorem, we conclude that the variational solution $u(t,x)$ is generated by a cohomology class of degree given by the Morse index of the quadratic form $\frac{1}{2}\langle Ap,p \rangle$ involved in (\ref{HAM}). 

\subsection{The convex case} \label{convex}
This section is devoted to detect the Markovian property for $p$-convex Hamiltonians. We emphasize that the proof is given by solely using minmax techniques.
\begin{Prop} \label{proof} Let $H(t,x,p) \in C^2([0,T] \times T^*\mathbb{T}^k)$ be $p$-convex. The corresponding $C^0$-variational solution is Markovian:
$$J^{t_3,t_2} \circ J^{t_2,t_1}(\sigma) = J^{t_3,t_1}(\sigma),$$
$\forall 0 \le t_1 \le t_2 \le t_3 \le T$ and $\forall \sigma \in C^{0}(\mathbb{T}^k)$.
\end{Prop}
\textit{Proof.} Through the proof, we indicate by $S_{s}^t$ the GFQI for the Hamiltonian flow $\phi^{t,s}_H$. For an arbitrary $C^1$ regularizing sequence $\sigma_n \stackrel{C^0}{\to} \sigma$, let (see (\ref{ESSE t}) and (\ref{min})):
\begin{equation} \label{alfa n}
\alpha_n (t,x) = \min_{(\xi,U)} [\sigma_n(\xi) + S_{t_1}^t(\xi,x;U)], \qquad t \in [t_1,T]
\end{equation}
and indicate by $\alpha \in C^0([t_1,T] \times N)$ the corresponding $C^0$-variational solution on $[t_1,T]$:
\begin{equation} \label{Uno}
\lim_{n \to +\infty} \| \alpha_n -\alpha \|_{C^0} = 0 \qquad \text{on } [t_1,T] \times N,
\end{equation}
so that $J^{t_2,t_1}(\sigma)(\cdot) = \alpha(t_2,\cdot)$. \\
In order to construct $J^{t_3,t_2} \circ J^{t_2,t_1}(\sigma)(\cdot)$, let us denote by 
\begin{equation} \label{APPROX-DATO}
f_n(\cdot) \stackrel{C^0}{\to} \alpha(t_2,\cdot)
\end{equation}
 an arbitrary approximating sequence, and
\begin{equation} \label{gamma n}
\gamma_n (t,x) = \min_{(\xi_1,U_1)} [f_n(\xi_1) + S_{t_2}^t(\xi_1,x;U_1)], \qquad t \in [t_2,T].
\end{equation}
Therefore, the variational solution $\gamma \in C^0([t_2,T] \times N)$ to the Cauchy problem starting from $\gamma(t_2, \cdot) = \alpha(t_2,\cdot)$ is given by:
$$
\lim_{n \to +\infty} \| \gamma_n - \gamma \|_{C^0} = 0 \qquad \text{on } [t_2,T] \times N,
$$
and $\gamma(t_3,\cdot) = J^{t_3,t_2} \circ J^{t_2,t_1} (\sigma)(\cdot)$. \\
Taking into account  (\ref{gamma n}) and (\ref{alfa n}), we detail now $\gamma_n(t,x)$:
$$\gamma_n(t,x) = \min_{(\xi_1,U_1)} \left[
f_n(\xi_1) + S_{t_2}^t(\xi_1,x;U_1)
\right]$$
$$= \min_{(\xi_1,U_1)} 
\left[ f_n(\xi_1) + S_{t_2}^t (\xi_1,x;U_1) - \alpha_n(t_2,\xi_1) + \alpha_n(t_2,\xi_1) \right]$$
$$= \min_{(\xi_1,U_1)} \left[f_n\left(\xi_1\right) + S_{t_2}^t\left(\xi_1,x;U_1\right) - \alpha_n\left(t_2,\xi_1\right) + \min_{(\xi,U)} \left[\sigma_n(\xi) + S_{t_1}^{t_2}(\xi,\xi_1;U)\right]\right]$$
\begin{equation} \label{FFF}
= \min_{(\xi_1,\xi,U_1,U)} \left[ \sigma_n(\xi) +
f_n\left(\xi_1\right) - \alpha_n\left(t_2,\xi_1\right) + S_{t_1}^{t_2}\left(\xi, \xi_1; U\right) + S_{t_2}^t\left(\xi_1,x;U_1\right) \right]
\end{equation}
$$= \min_{(\xi_1,\xi,U_1,U)} \left[ \sigma_n(\xi) + (S_n)_{t_1}^{t})(\xi,x;\xi_1,U_1,U) \right],$$
where $(S_n)_{t_1}^t$ involves all the above terms but $\sigma_n$. We claim that: 
$$\lim_{n \to +\infty} \| (S_n)_{t_1}^{t} - S_{t_1}^t \|_{C^0} = 0 \qquad \text{on } [t_2, T] \times N.$$
In fact, in view of the composition rule for generating functions (cfr. Lemma \ref{RULE} in the Appendix), 
$S_{t_1}^{t_2}\left(\xi, \xi_1; U\right) + S_{t_2}^t\left(\xi_1,x;U_1\right)$ in (\ref{FFF}) generates $\phi_H^{t,t_1}$ for $t \ge t_2$, and the convergence follows from (\ref{APPROX-DATO}). Consequently, see also Proposition III in \cite{VIO}:
\begin{equation} \label{GAMMA_N}
\lim_{n \to +\infty} \| \gamma_n - \alpha_n \|_{C^0} = 0 \qquad \text{on } [t_2,T] \times N.
\end{equation}
Therefore, from 
$$\| \gamma_n - \alpha \|_{C^0} \le \| \gamma_n - \alpha_n \|_{C^0} +  \| \alpha_n - \alpha \|_{C^0}$$
and by using (\ref{GAMMA_N}) and (\ref{Uno}), we achieve that $\alpha = \gamma$ on $[t_2, T]$ and for $t = t_3$ we obtain the thesis. \hfill $\Box$ \\ \\
The proof of the previous proposition prompts generalizations towards some non-convex cases, which will be discussed here below.

\subsection{The convex-concave case} \label{ex}
Let $j$ be an integer, $0 \le j \le k$, and for any $(x,p) \in T^*\mathbb{T}^k$ set:
$$x=(x_1,x_2) \in \mathbb{T}^j \times \mathbb{T}^{k-j}, \qquad p=(p_1,p_2) \in \mathbb{R}^j \times \mathbb{R}^{k-j}$$
and assume:
\begin{equation} \label{ham}
H(t,x,p) = H_1(t,x_1,p_1) + H_2(t,x_2,p_2) \in C^{2}([0,T] \times T^*\mathbb{T}^k),
\end{equation}
where $H_1(t,x_1,p_1)$ is $p$-convex and $H_2(t,x_2,p_2)$ is $p$-concave. In this section, we denote by $(S_1)^t_0$ and $(S_2)^t_0$ the generating functions for $\phi^{t,0}_{H_1}$ and $\phi^{t,0}_{H_2}$ respectively. \\
\indent We first suppose that the initial datum $\sigma \in C^1(\mathbb{T}^k)$ is of the form:
\begin{equation} \label{sig}
\sigma(x) = \sigma_1(x_1) + \sigma_2(x_2),
\end{equation} 
and we prove the following
\begin{Prop} \label{inf-sup} The variational solution $u(t,x)$ to the Hamilton-Jacobi equation with Hamiltonian (\ref{ham}) and initial datum (\ref{sig}) is given for all $(t,x) \in [0,T] \times \mathbb{T}^k$ by:
\begin{equation} \label{sol1}
u(t,x) =  \min_{(\xi_1,U_1)}\left[ \sigma_1 (\xi_1) + (S_1)^t_0 (\xi_1,x_1;U_1) \right] + \max_{(\xi_2,U_2)} \left[ \sigma_2(\xi_2) + (S_2)^t_0(\xi_2,x_2;U_2) \right].
\end{equation}
\end{Prop}
\textit{Proof.} The Lagrangian wavefronts $\phi^{t,0}_{H_1}(\text{Im}(d \sigma_1))$ and $\phi^{t,0}_{H_2}(\text{Im}(d \sigma_2))$ admit the GFQI:
$${\bf{{\cal{S}}}}_1(t,x_1;\xi_1,U_1) := \sigma_1(\xi_1) + (S_1)^t_0(\xi_1,x_1;U_1)$$ 
and 
$${\bf{{\cal{S}}}}_2(t,x_2;\xi_2,U_2) := \sigma_2(\xi_2) + (S_2)^t_0(\xi_2,x_2;U_2),$$
with positive and negative defined quadratic form respectively (see Theorem \ref{AA}). Moreover, since $H$ and $\sigma$ separate the variables $x_1$ and $x_2$, $\phi^{t,0}_{H}(\text{Im}(d \sigma))$ is generated by the GFQI with parameters $(\xi,U) = (\xi_1,\xi_2,U_1,U_2)$:
\begin{equation} \label{ESSE-12}
{\bf{{\cal{S}}}}(t,x;\xi,U) := {\bf{{\cal{S}}}}_1(t,x_1;\xi_1,U_1) + {\bf{{\cal{S}}}}_2(t,x_2;\xi_2,U_2)
\end{equation}
and related quadratic form $Q(\xi,U) = (\xi_1^2 + U_1^2) - (\xi_2^2 + U_2^2)$. \\
Let $c > 0$ be large enough so that the corresponding sublevel sets ${\bf{{\cal{S}}}}^{\pm c}$ are respectively connected and unconnected. Keeping in mind that the minmax critical value $u(t,x)$ marks a metamorphosis of the topology of the sublevel sets for $\bf{{\cal{S}}}$, we proceed by giving upper and lower bounds for $u(t,x)$. \\
From the one hand, by decreasing $\lambda \le c$, as long as the sublevel set ${\bf{{\cal{S}}}}^{\lambda}$, given by
\begin{equation} \label{in}
\sigma_1(\xi_1) + (S_1)^t_0(\xi_1,x_1;U_1) + \sigma_2(\xi_2) + (S_2)^t_0(\xi_2,x_2;U_2) \le \lambda
\end{equation}
remains connected, for all $(\xi_2,U_2)$ there exists at least one $(\xi_1,U_1)$ such that the inequality (\ref{in}) is satisfied. In particular, there exists $(\xi_1,U_1)$ such that
$$\sigma_1 (\xi_1) + (S_1)^t_0(\xi_1,x_1;U_1) + \max_{(\xi_2,U_2)} \left[ \sigma_2(\xi_2) + (S_2)^t_0(\xi_2,x_2;U_2) \right] \le \lambda.$$
Therefore:
\begin{equation} \label{prima}
u(t,x) \le \min_{(\xi_1,U_1)}\left[ \sigma_1 (\xi_1) + (S_1)^t_0(\xi_1,x_1;U_1) \right] + \max_{(\xi_2,U_2)} \left[ \sigma_2(\xi_2) + (S_2)^t_0(\xi_2,x_2;U_2) \right].
\end{equation}
From the other hand, by increasing $\lambda \ge -c$, as long as the sublevel set ${\bf{{\cal{S}}}}^{\lambda}$ is unconnected, there exists $(\xi_2,U_2)$ such that for all $(\xi_1,U_1)$ the previous inequality (\ref{in}) is unsatisfied. In particular, it exists $(\xi_2,U_2)$ such that
$$\min_{(\xi_1,U_1)} \left[ \sigma_1 (\xi_1) + (S_1)^t_0(\xi_1,x_1;U_1) \right] + \sigma_2(\xi_2) + (S_2)^t_0(\xi_2,x_2;U_2) > \lambda,$$
and consequently:
\begin{equation} \label{seconda}
u(t,x) \ge \min_{(\xi_1,U_1)}\left[ \sigma_1 (\xi_1) + (S_1)^t_0(\xi_1,x_1;U_1) \right] + \max_{(\xi_2,U_2)} \left[ \sigma_2(\xi_2) + (S_2)^t_0(\xi_2,x_2;U_2) \right].
\end{equation}
From the estimates (\ref{prima}) and (\ref{seconda}), we obtain the thesis. \hfill $\Box$ \\ \\
As a consequence of the explicit formula (\ref{sol1}), the variational solution to the Hamilton-Jacobi equation with convex-concave Hamiltonian (\ref{ham}) and initial datum (\ref{sig}) results Markovian. In fact, since (\ref{sol1}) is the superimposition of the solutions to the Cauchy problems:
$$\begin{cases}
\frac{\partial u}{\partial t}(t,x) + H_i(t,x_i,\frac{\partial u}{\partial x}(t,x_i)) = 0 \\ \\
u(0,x_i) = \sigma_i(x_i)
\end{cases}$$
$i = 1,2$, the arguments of the proof of Proposition \ref{proof} still hold. \\
\indent We conclude this section by considering a generic initial datum, that is not split up in the form (\ref{sig}). The Lagrangian wavefront $\phi^{t,0}_H(\text{Im}(d\sigma))$ is now generated by:
\begin{equation} \label{ESSE-MISTA}
{\bf{{\cal{S}}}}(t,x;\xi,U) :=  \sigma(\xi_1,\xi_2) + (S_1)^t_0(\xi_1,x_1;U_1) + (S_2)^t_0(\xi_2,x_2;U_2).
\end{equation}
In such a case, the same arguments on the topology of the sublevel sets used in the proof of Proposition \ref{inf-sup} lead to the following explicit pointwise upper and lower bounds for the corresponding variational solution.
\begin{Prop} \label{ULTI} The variational solution $u(t,x)$ to the Hamilton-Jacobi equation with Hamiltonian (\ref{ham}) and continuous initial datum $\sigma \in C^1(\mathbb{T}^k)$ satisfy for all $(t,x) \in [0,T] \times \mathbb{T}^k$
\begin{equation} \label{sol11}
\max_{(\xi_2,U_2)}  \min_{(\xi_1,U_1)} \left[ {\bf{{\cal{S}}}}(t,x;\xi,U) \right] \le u(t,x) \le \min_{(\xi_1,U_1)} \max_{(\xi_2,U_2)} \left[ {\bf{{\cal{S}}}}(t,x;\xi,U) \right],
\end{equation}
where ${\bf{{\cal{S}}}}(t,x;\xi,U)$ is given by (\ref{ESSE-MISTA}).
\end{Prop}
The previous proposition restores in a genuine minmax framework a rather general Hopf-type estimate and gives a representation formula for the variational solution whenever the first and the last terms are equal. This trivially occurs for $j = k$ or $j = 0$, in such a case $u(t,x)$ reduces to the formulas (\ref{min}) and (\ref{max}) respectively. A more interesting case has been just discussed in Proposition \ref{inf-sup}. \\ \\
\textbf{Remark} Similar Hopf-type estimates hold for viscosity solutions, see \cite{BF} and \cite{BO}. In particular, the paper \cite{BF} treats the case of a strictly integrable convex-concave Hamiltonian $H(p) = H_1(p_1) + H_2(p_2)$ and presents (\ref{sol1}) and (\ref{sol11}) for viscosity solutions. Therefore, the Proposition \ref{inf-sup} establishes definitively the coincidence for variational and viscous solutions to Hamilton-Jacobi problems for convex-concave Hamiltonians $H(p)$ and initial data separating the variables. However, we stress that the present estimates (\ref{sol11}) for variational solutions do hold for a larger class of Hamiltonians, possibly depending on $(t,x)$. We finally note that the papers \cite{BF} and \cite{BO} also present upper and lower bounds for viscosity solutions in the case where the convex-concave assumptions are not on the Hamiltonian function but on the initial datum: this is consistent with the viscosity solutions theory, which does not operate necessarily with compactly supported initial data and GFQI.

\section{Appendix: Construction of global GFQI} \label{otto}
For the large class of Hamiltonians $H(t,x,p) = \frac{1}{2}\langle Ap,p \rangle + V(t,x,p)$ introduced  in Section \ref{sette}, here we explicitly construct the finite parameters GFQI globally describing the geometric solution $L$. \\
\indent Denoting by $H_t(x,p) = H(t,x,p)$, we lift $H_t(x,p)$ to the covering space $[0,T] \times \mathbb{R}^{2k}$ and indicate the components of its flow by  $(x_s^t,p_s^t) := \phi_H^{t,s}(X,P)$. The following proposition is a consequence of the above assumptions $(\star)$ on $H_t(x,p)$ (see also Proposition 2.3.3 in \cite{CH2}). 
\begin{Prop} It exists $\varepsilon > 0$ small enough so that every application: 
\begin{equation} \label{ACCA}
h_{s}^{t}: (X,P) \mapsto (X,x_s^t), \qquad 0 < |t-s| < \varepsilon
\end{equation}
is a Lipschitz diffeomorphism of $\mathbb{R}^{2k}$. 
\end{Prop}
\textit{Proof.} Note that (\ref{ACCA}) is degenerate for $s = t$. To overcome this fact, we use Chaperon's trick by introducing the following isomorphism:
$$L_{s}^{t}(X,x) := \left( X, \frac{x - X}{t - s} \right),$$
and we equivalently prove the assertion for $k_{s}^{t} := L_{s}^{t} \circ h_{s}^{t}$. The application $k_{s}^{t}$ is well-defined also for $t = s$. In fact, Hamilton's equations gives:
\begin{equation} \label{kappa}
k_{s}^{t}(X,P) = \left( X, \frac{1}{t - s} \int^{t}_{s} D_pH_{\tau}(\phi^{\tau}_{s}(X,P)) d\tau \right),
\end{equation}
and therefore $k_{s}^{s}(X,P) = (X,D_pH_{s}(X,P))$, which is, precisely, the Legendre transformation. Conditions on $H_t(x,p)$ guarantee us that $k^s_s$ is a Lipschitz diffeomorphism for every $s \in [0,T]$ and that the Lipschitz constants of $k^s_s$ and $(k_s^s)^{-1}$ are bounded independently of $s$. We estimate now the Lipschitz constant of $k^t_s - k^s_s$. We have:
$$(k_s^t - k_s^s)(X+x,P+p) - (k_s^t - k_s^s)(X,P) =$$
$$= \left(0,\frac{1}{t-s}\int_s^t\left[D_pH_{\tau}(\phi_s^{\tau}(X+x,P+p)) - D_pH_{\tau}(X+x,P+p)\right] d\tau\right) -$$
$$- \left(0,\frac{1}{t-s}\int_s^t\left[D_pH_{\tau}(\phi_s^{\tau}(X,P)) - D_pH_{\tau}(X,P)\right] d\tau\right) \to (0,0)$$
for $t \to s$. Moreover: 
$$\left\| \left( D_p H_{\tau} \circ \phi^\tau_s - D_p H_{\tau} \right)(X+x,P+p) - \left( D_p H_{\tau} \circ \phi^\tau_s - D_p H_{\tau} \right)(X,P) \right\|_{\mathbb{R}^{k}} \le$$
$$\le \text{Lip}(D_p H_{\tau})\left\| \phi^{\tau}_s(X+x,P+p) - \phi^{\tau}_s(X,P) \right\|_{\mathbb{R}^{2k}} + \text{Lip}(D_pH_{\tau})\left\| (x,p)\right\|_{\mathbb{R}^{2k}} \le$$
$$\le \left( \text{Lip}\left(D_pH_{\tau}\right)\text{Lip}(\phi^{\tau}_s) + \text{Lip}\left(D_pH_s\right) \right) \left\| (x,p) \right\|_{\mathbb{R}^{2k}},$$
where, in our hypothesis, the Lipschitz constants appearing in the last member are uniformly bounded with respect to $s$ and $\tau$ in $[0,T]$. As a consequence, it results that $\text{Lip}(k^t_s - k^s_s) \to 0$ for $t \to s$, \textit{uniformly} with respect to $s$. Hence, it exists $\varepsilon > 0$ such that
$$\text{Lip}(\text{Id} - k^t_s \circ (k_s^s)^{-1}) = \text{Lip}((k^s_s - k^t_s) \circ (k_s^s)^{-1}) \le$$
$$\le \text{Lip}(k^s_s - k^t_s) \sup_{s \in [0,T]} \text{Lip}((k_s^s)^{-1}) < 1$$
for $0 < |t-s| < \varepsilon$. It is now well-known that contractible perturbations of the identity (that is, $\text{Id} - u$ where $\text{Lip}(u) < 1$) are Lipschitz diffeomorphisms. Applying this result to $u = (\text{Id} - k^t_s \circ (k_s^s)^{-1})$, we obtain that the application $\text{Id} - u = k^t_s \circ (k_s^s)^{-1}$ is a Lipschitz diffeomorphism; and the same is true for $k^t_s$. \hfill $\Box$ \\ \\
As a consequence, quoting for example \cite{FM} and \cite{GO}, the flow:
$$\phi^{t,s}_H: (X,P) \to (x_s^t,p_s^t),$$
$0 < |t-s| < \varepsilon$, provides a symplectic \textit{twist} diffeomorphism. In the sequel, we alternatively use the expression symplectic diffeomorphism and canonical transformation. \\
\indent Fixed $t \in [0,T]$, we divide the interval $[0,t]$ in an appropriate number $N+1$ of sub-intervals $[t_j,t_{j+1}]$, $0 \le j \le N$, of length $\varepsilon = \frac{t}{N+1}$:
$$0=t_0 < \ldots < t_{N+1} = t, \qquad N \in \mathbb{N}$$ 
and we define
\begin{equation} \label{notazioni}
(X_{j+1},P_{j+1}) := \phi^{t_{j+1},t_j}_H(X_j,P_j).
\end{equation}
As a consequence of the previous proposition, it exists $N$ large enough such that \textit{every} application: $(X_j,P_j) \mapsto (X_j,X_{j+1})$ is a global diffeomorphism of $\mathbb{R}^{2k}$. Hence, one can obtain the explicit expression for $P_j$: 
$$P_j = \tilde{P}_j(X_j,X_{j+1}),$$
and accordingly use the variables $(X_j,X_{j+1})$ in order to describe the canonical transformation given by the flow $\phi^{t_{j+1},t_j}_H$. The fact that the map $\phi^{t_{j+1},t_j}_H$ is symplectic and the triviality of the deRham cohomology group $H^1(\mathbb{R}^{2k}) = \{ 0 \}$ imply the existence of a generating function $S^{t_{j+1}}_{t_j}(X_j,X_{j+1})$ such that:
\begin{equation} \label{REL}
\begin{cases}
P_j = - \frac{\partial}{\partial X_j} S^{t_{j+1}}_{t_j}(X_j,X_{j+1}) \\ \\
P_{j+1} = \frac{\partial}{\partial X_{j+1}} S^{t_{j+1}}_{t_j}(X_j,X_{j+1})
\end{cases}
\end{equation}
Moreover, $S^{t_{j+1}}_{t_j}$ comes essentially from the Hamilton-Helmoltz functional: 
$$S_{t_j}^{t_{j+1}}(X_j,X_{j+1}) =$$
\begin{equation} \label{EsseJ} 
= \left[ \int^{t_{j+1}}_{t_j} \left( p^\tau_{t_j}\left(X_j,P_j\right) \frac{d}{d \tau} x^\tau_{t_j}\left(X_j,P_j\right) - H_{\tau}\left(\phi^{\tau,t_j}_H\left(X_j,P_j\right)\right) \right) d\tau \right]_{\arrowvert_{P_j = \tilde{P}_j(X_j,X_{j+1})}}
\end{equation} 
see \cite{CH2} for a detailed proof of this fact. \\
\indent The Hamilton principle function relative to the canonical transformation given by the flow $\phi^{t,0}_H$ can now be easily computed by using the next version of a popular composition rule for generating functions of canonical transformations (see for example \cite{LL1}, \cite{LL2} and \cite{Be}). 
\begin{Lem} \label{RULE} Let
$$\phi: (Q,P) \mapsto (q,p), \qquad \psi: (X,Y) \mapsto (x,y)$$
be two canonical transformations generated respectively by $S(Q,q)$ and $F(X,x)$. Then the canonical transformation 
$$\psi \circ \phi: (Q,P) \mapsto (x,y)$$ 
is generated by
\begin{equation} \label{CR}
G(Q,x;w) := S(Q,w) + F(w,x).
\end{equation} 
\end{Lem} 
\textit{Proof.} It is sufficient to write down the relations for $G(Q,x;w)$: $P = -\frac{\partial G}{\partial Q}$, $y = \frac{\partial G}{\partial x}$, $0 = \frac{\partial G}{\partial w}$. \hfill $\Box$ \\
\begin{Prop} \label{teoremaprincipale} Let $(X,x) := (X_0,X_{N+1})$ and $U := (X_{j})_{1 \le j \le N} \in \mathbb{R}^{Nk}$.  The canonical transformation given by the flow 
$$\phi^{t,0}_H: (X,P) \mapsto (x,p)$$
is generated by
\begin{equation} \label{FG}
S^{t}_0(X,x;U) := \sum_{j=0}^{N}S^{t_{j+1}}_{t_j}(X_j,X_{j+1}).
\end{equation}
In other words, the graph of $\phi^{t,0}_H$ is the set of points:
\begin{equation} \label{Set}
\left\{
\left(
\left(X, Y = - \frac{\partial S^{t}_0}{\partial X}\right),
\left(x,p = \frac{\partial S^{t}_0}{\partial x}\right)
\right): \ \frac{\partial S^{t}_0}{\partial U} = 0
\right\} \subseteq \mathbb{R}^{2k} \times \mathbb{R}^{2k}.
\end{equation}
\end{Prop}
Theorem \ref{AA} in Section \ref{sette} is now a straightforward consequence of the previous proposition. \\ \\
\textit{Proof of Theorem \ref{AA}.} A simple direct computation --in particular, the stationarization with respect to $\xi$ connecting the impulses of $\text{Im}(d\sigma)$ with the starting impulses of $\phi^{t,0}_H$-- proves that the function $S^t(x;\xi,U) := \sigma(\xi) + S^t_0(\xi,x;U)$ generates the Lagrangian submanifold $\phi^{t,0}_H(\text{Im}(d\sigma))$. We check below the quadraticity at infinity property of $S^t(x;\xi,U)$ with respect to the parameters $(\xi,U)$: this is a crucial step in order to catch the minmax critical point in the Lusternik-Schnirelman format. \\
Using the notation $(X_{j+1},P_{j+1}) = \phi^{t_{j+1},t_j}_H(X_j,P_j)$, we refer to a partition of the interval $[0,t]$, $0=t_0 < \ldots < t_{N+1} = t$ and $t_{j+1} - t_j = \varepsilon$, so that every application $(X_j,P_j) \mapsto (X_j,X_{j+1})$ is a global diffeomorphism of $\mathbb{R}^{2k}$ and (\ref{REL}) holds. Since $V(t,x,p)$ is compactly supported in the $p$ variables, we have that $H(t,x,p) = \frac{1}{2}\langle Ap,p \rangle$ for $p$ outside a compact set. As a consequence, we can assume that, up to the three main operations described in Section \ref{Aggiunta}, every generating function $S_{t_j}^{t_{j+1}}(X_j,X_{j+1})$ coincides, for large $\| \frac{A^{-1}(X_{j+1} - X_{j})}{\varepsilon} \|_{\mathbb{R}^k}$, to 
$$S_{t_j}^{t_{j+1}}(X_j,X_{j+1}) = \frac{1}{2 \varepsilon} \Big\langle A^{-1}(X_{j+1} - X_{j}), (X_{j+1} - X_{j}) \Big\rangle.$$
In fact, we obtain:
$$\begin{cases}
P_j = -\frac{\partial}{\partial X_j} S_{t_j}^{t_{j+1}}(X_j,X_{j+1}) = \frac{1}{\varepsilon} A^{-1}(X_{j+1} - X_j) \\ \\
P_{j+1} = \frac{\partial}{\partial X_{j+1}} S_{t_j}^{t_{j+1}}(X_j,X_{j+1}) = \frac{1}{\varepsilon} A^{-1}(X_{j+1} - X_j)
\end{cases}$$
that is, $P_{j+1} = P_j = \frac{1}{\varepsilon}A^{-1}(X_{j+1} - X_{j})$ and $X_{j+1} = X_j + \varepsilon A P_j$. We apply now proposition \ref{teoremaprincipale}: denoting $(\xi,x) := (X_0,X_{N+1})$ and $U = (X_j)_{1 \le j \le N}$, we have that
$$S^t(x;\xi,U) = \sigma(\xi) + S_0^t(\xi,x;U) =  \sigma(X_0) + \sum_{j=0}^N S_{t_j}^{t_{j+1}}(X_j,X_{j+1}).$$
Finally, by introducing the following change of variables:
\begin{equation*} \label{cambio}
\xi_j := A^{-1}(X_{j+1} - X_{j}) \qquad \text{for} \ 0 \le j \le N,
\end{equation*}
the resulting function: 
\begin{equation}
S^t(x; (\xi_j)_{0 \le j \le N}) = \sigma\left(x - \sum_{j=0}^N A \xi_j\right) + \frac{1}{2 \varepsilon} \sum_{j = 0}^N \langle A \xi_j, \xi_j \rangle
\end{equation}
is $\mathbb{Z}^k$-periodic with respect to $x$ and is quadratic at infinity, with nondegenerate quadratic form given by:
\begin{equation}
\frac{1}{2\varepsilon} \sum_{j = 0}^N \langle A \xi_j, \xi_j \rangle = \frac{1}{2 \varepsilon}(\xi_1, \ldots, \xi_N)\textbf{A}(\xi_1, \ldots, \xi_N)^t.
\end{equation} \hfill $\Box$ \\ \\
\textbf{Acknowledgments} O. Bernardi has been supported by the project ``Te\-cni\-che va\-ria\-zio\-na\-li e PDE in to\-po\-lo\-gia sim\-plet\-ti\-ca appli\-ca\-zio\-ni fi\-si\-co-ma\-te\-ma\-ti\-che'' of Gruppo Nazionale di Fisica Matematica G.N.F.M.


\end{document}